# Modeling of Phenomena and Dynamic Logic of Phenomena

Boris Kovalerchuk[1], Leonid Perlovsky[2], Gregory Wheeler[3]

[1]Central Washington University, Dept. of Computer Science, USA, borisk@cwu.du
[2]Harvard University and the Air Force Research Laboratory, USA, leonid@seas.harvard.edu
[3]New University of Lisbon, CENTRIA – AI Center, Dept. of Computer Science, Portugal, grw@fct.unl.pt

*Abstract*—Modeling of complex phenomena such as the mind presents tremendous computational complexity challenges. Modeling field theory (MFT) addresses these challenges in a non-traditional way. The main idea behind MFT is to match levels of uncertainty of the model (also, problem or theory) with levels of uncertainty of the evaluation criterion used to identify that model. When a model becomes more certain, then the evaluation criterion is adjusted dynamically to match that change to the model. This process is called the Dynamic Logic of Phenomena (DLP) for model construction *and it mimics processes of the mind and natural evolution. This paper provides a formal description of DLP by specifying its syntax, semantics, and reasoning system. We also outline links between DLP and other logical approaches. Computational complexity issues that motivate this work are presented using an example of polynomial models.*

## I. Introduction

There are two current trends within the modeling of physical phenomena and within logic. In modeling physical phenomena, the trend is to put in additional logical structures to classical mathematical techniques, whereas in logic it is to add "dynamics" to the logic [12-17], with an aim to represent and reason about *actions* rather than static *propositions.* The subjects in this area include Action Logic, Arrow Logic, Game Logic, Semantic games, Dialogue Logic, Belief Revision, Dynamic Epistemic Logic, Hoare logic, Dynamic Logic, Linear Logic, labeled transition systems, Petri nets, Process Algebra, Automata Theory, Game Semantics, coalgebras, among others. [13].

These two complimentary trends can be very beneficial for both areas. In [13] "dynamification" of logic is developed to model quantum phenomena, and this paper develops the dynamic approach for modeling other problems where the computational complexity of finding the solution is a critical issue.

While at the global level these complimentary trends exist, for one approach to benefit the other, both need to be close enough to each other in specific tasks and goals.
We start by discussing the differences between dynamic epistemic logic [16, 19] and DLP. H. Leitgeb has pointed out that the aim of dynamic epistemic logic and the like is to put logical operators with a dynamic interpretation *into* one's formal object language. The aim of DLP, however, is to specify a particular dynamics of learning and related concepts *in the meta-language.* So, since DLP locates the dynamics in the meta-language and dynamic epistemic logic locates them in the object language, the logical resources of dynamic epistemic logic are not yet of help to DLP [18].

To illustrate this difference, consider an example from [24]: "A man loves Annie. He is rich." Two interpretations are possible:
(1) $(\exists x)$(man x & x loves Annie) & x is rich.
(2) $(\exists x)$(man x & x loves Annie & x is rich).
In the first sentence, the existential quantifier is applied only to the first sentence, but in the second sentence, it ranges over both. The two interpretations of "A man loves Annie. He is rich" show that the existential quantifier and operations can be dynamic, that is, they can have different interpretations within the object language. According to J. van Benthem [20], DLP follows dynamic systems tradition. This tradition can be traced to Ernst Max who viewed organisms as dynamic systems that have innate tendencies to self-regulation and equilibrium. When equilibrium is disturbed, which can happen on a variety of levels, the organism works to form a new equilibrium [21]. Note that the common tools to model such phenomena are differential equations rather than logic. A related view of dynamic systems comes from studies in Computational Theory of Mind (CMT) [22]. According to [23], cognitive processes are not rule-governed sequences of discrete symbolic states, but continuously evolving total states of dynamic systems determined by continuous, simultaneous and mutually determining states of the systems' components (i.e., state variables or parameters). In section 11, we outline a way to shrink the gap between these two basic approaches.

To provide a formal description of Dynamic Logic of Phenomena (DLP), we start by comparing the background definitions of DLP to logic model theory. Section 2 establishes concepts of uncertainty, generality, and simplicity for models, and defines evaluation criteria. Section 3 defines a partial order of models. Section 4 provides examples of uncertainty and generality of polynomial models. Section 5 formalizes similarity maximization. Section 6 defines DLP parameterization using the theory of monotone Boolean functions. Section 7 defines the search process. Section 8 presents how DLP processes can be visualized. Section 9 provides a formal description of DLP. Sections 10 and 11 outline the links between DLP and other dynamic logics. Section 12 summarizes the paper and discusses future research. The Appendix describes computation complexity issues that motivate the paper.

We start by defining the concept of *empirical data* relevant to modeling field theory (MFT) [8], [9], and supply an interpretation in logical terms.

**Empirical data, E** in MFT is any data to identify a model.

In logical terms, we define empirical data as a pair, $E=\langle A, \Omega\rangle$, where A is a set of objects, and $\Omega=\{P_i\}$ is a set of predicates $P_i$ of arity $n_i$, e.g., $P_1(x,y)$ means that length of x is no less than the length of y, $l(x) \geq l(y)$.

<u>Definition</u>. A pair $<A, \Omega>$ is called an **empirical system** [6].

<u>Definition.</u> A pair $<A, \Omega>$ is called **a *model*** in logic [7]. Often it is considered as a model of some system of axioms T.

Tarski proposed the name 'model theory' in 1954, although a variety of other names are also used, including *relational system* [6], and a *protocol of the experiment*. We call a Tarskian models **logic models** or **Lmodels** [10] to distinguish them from models in MFT.

The concept of a model, M, in MFT concerns a model of reality, which we will call a **model of phenomena** or **Pmodel**. In logical formalization, a Pmodel can be matched with an **axiom system** T.

<u>Definition</u>. A **system of axioms** T is a set of closed formulas (sentences) in the signature of the underlying language, e.g.,

$$\forall x_i \exists x_j P_1(x_i, x_j).$$

The concept of model is treated very differently in MFT than it is in mathematical logic. Logic may be thought to go from a *very formal* (syntactical) axiomatic system T to a more real or concrete model $\mathbf{A_T} = <A_T, \Omega>$ of that formal system T. MFT goes in the other direction, from a *very informal* reality to more formal models. As a result, the concepts of model are quite different in the two theories. Empirical data in MFT is a model $\mathbf{E}=<A, \Omega>$ in mathematical logic, if we interpret empirical data as an empirical system E [6]. On the other hand, the model of phenomena is *not* a model in logic; instead, Pmodels are akin to *a set of axioms* about the class of logic models. This type of difference was well described in [10]:

> "To model a *phenomenon* is to construct a formal theory that describes and explains it. In a closely related sense, you *model* a system or structure that you plan to build, by writing a description of it. These are very different senses of 'model' from that in model theory: the '**model' of the phenomenon** or the system is not a structure but a *theory*, often in a formal language."

Thus, we will use terms that have been already introduced above: **Pmodel** for a model of phenomenon and **Lmodel** for a logic model.

The next MFT concept is a **similarity (or correspondence) measure**, L(M,E), between empirical data E and an a-priori model M that is assigned individually to each pair (M,E):

$$L: \{(M,E)\} \to R,$$

where R is a set or real numbers. In logic the closest to it is a statement that $<A,\Omega>$ is a model of the system of the axioms T.

In logical terms, L maps a theory M and Lmodel E to R.

<u>Definition</u>. Pair $E = <A,\Omega>$ is an **Lmodel of the system of the axioms** T if every formula from T is true on E.

<u>Definition</u>. **Boolean similarity measure** B(T,E) is defined to be equal to 1, B(T,E)=1, If M is an Lmodel of T, else B(T,E)=0.

## II. SEMANTIC CONCEPTS OF UNCERTAINTY, GENERALITY AND SIMPLICITY

### A. Uncertainty, generality, and simplicity relations between P-models

Below we introduce the concepts of uncertainty, generality, and simplicity relations. These concepts can be specified for both logic and MFT models.

An **uncertainty relation between Pmodels** is denoted by $\geq_{Mu}$, and the sentence $M_i \geq_{Mu} M_j$ is read: "Model $M_i$ is equal in uncertainty or *more uncertain* than model $M_j$". In other words, model $M_j$ is equal in certainty or more certain than model $M_i$, and model $M_j$ is no less certain than model $M_i$. This relation is a partial order. If $M_i >_{Mu} M_j$ then we simply say that $M_j$ is more certain than $M_i$.

A **generality relation between Pmodels** is denoted by $\geq_{Mg}$ and relation $M_i \geq_{Mg} M_j$ is read: "Model $M_j$ is a *specialization* of the model $M_i$" or "Model $M_i$ is a *generalization* of the model $M_j$". This relation also is a partial order.

A **simplicity relation between Pmodel** is denoted by $\geq_{Ms}$ and relation $M_i \geq_{Ms} M_j$ is read: "Model $M_i$ is equal in simplicity of *simpler* than Model $M_j$". This relation also is a partial order.

For Pmodels that are represented as a system of axioms, the generality relation can be defined as follows.

<u>Definition.</u> $T_i \geq_{gen} T_j$ if and only if $T_i \subseteq T_j$, i.e., system of axioms $T_i$ is equal to, or an extension of, the system of axioms $T_j$ if and only if every axiom in $T_i$ is an axiom in Tj.

### B. Uncertainty, generality and simplicity relations between similarity measures

An **uncertainty relation** between similarity measures is denoted by $\geq_{Lu}$, and $L_i \geq_{Lu} L_j$ is read: "Measure $L_i$ is equal in uncertainty or more uncertain than measure $L_j$". This is a partial order relation.

A **generality relation between similarity measures** is denoted by $\geq_{Lg}$, and $L_i \geq_{Lg} L_j$ is read either: "Measure $L_j$ is a *specialization* of measure $L_i$", or equivalently, "measure $L_i$ is a *generalization* of the measure $L_j$". This relation also is a partial order.

A **simplicity relation between similarity measures** is denoted by $\geq_{Ls}$, and relation $L_i \geq_{Ls} L_j$ is read: "Measure $L_j$ is equal in simplicity or *simpler* than measure $L_i$". This relation also is a partial order.

<u>Definition</u>. Mapping F between a **set of Pmodels** {M} and a set of similarity measures {L},

$$F: \{M\} \to \{L\},$$

is called a **match mapping** if F preserves uncertainty, generality, and simplicity relations between models and measures in the form of *homomorphism* from a relational system $< \{M\}, \geq_{Mg}, \geq_{Mu}, \geq_{Ms} >$ to a relational system

$< \{L\}, \geq_{Lg}, \geq_{Lu}, \geq_{Ls} >$, i.e.,

$\forall M_a, M_b \ (M_a \geq_{Mg} M_b \Rightarrow F(M_a) \geq_{Lg} F(M_b))$,

$\forall M_a, M_b \ (M_a \geq_{Mu} M_b \Rightarrow F(M_a) \geq_{Lu} F(M_b))$,

$\forall M_a, M_b \ (M_a \geq_{Ms} M_b \Rightarrow F(M_a) \geq_{Ls} F(M_b))$.

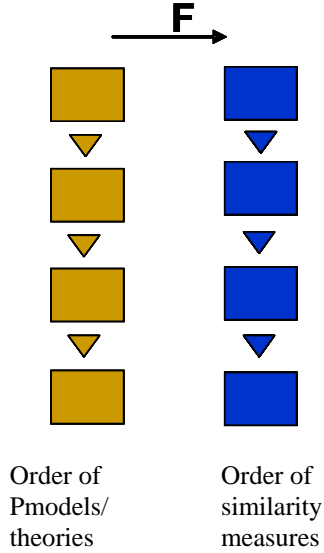

Figure 1. Mapping of Pmodels/logic theories and similarity measures

## II. PARTIAL ORDER OF PMODELS

Two different Pmodels can be at the same level of uncertainty ($M_1 =_u M_2$), one Pmodel can be more uncertain than another ($M_1 >_u M_2$), or Pmodels can be incomparable for uncertainty. We may define *model uncertainty* in such a way that two different models,

$$M_1: \forall x,y \ 2x^2 + 3y = 0$$

and

$$M_2: \forall x,y \ 5x + 4y^2 = 0$$

have the same level of uncertainty $M_1 =_u M_2$. The *number of unknown coefficients* is one of the possible ways to define the level of uncertainty. For $M_1$ and $M_2$, these numbers $m_1$ and $m_2$ are equal to zero. All coefficients are known and models are certain. In contrast model

$$M_3: \forall x,y \ 5x + by2 = 0.$$

has one unknown coefficient b and its measure of uncertainty $m_3$ is equal to 1.

Definition. *NUC measure of polynomial model uncertainty* is defined as the Number of Unknown Coefficients (NUC) in the model.

Consider Pmodel $M_3$. Based on the NUC measure, this Pmodel is *more uncertain* than Pmodel $M_2$, $M_3 >_{Mu} M_2$, because coefficient b in $M_3$ is unknown. NUC for $M_2$ is $n_2=0$ and NUC for $M_3$ is $n_3=1$ and $n_3>n_2$. In logical terms these Pmodels can be described as logic statements

$$M_1: \forall x,y \ P_1(x,y),$$

$$M_2: \forall x,y \ P_2(x,y),$$

$$M_3: \forall x,y \ P_3(x,y,b).$$

The *generality* relation between Pmodels $M_1$ and $M_2$ can also be defined. For instance, it can be the *highest power n* of the polynomial model. Both models $M_1$ and $M_2$ are quadratic with $n_1=n_2=2$ and, thus, both have the same generality.

Definition. *HP measure of polynomial model generality* is defined as the *Highest Power n* of the polynomial model.

Alternatively, we may look deeper and notice that $M_1$ contains $x^2$ and $M_2$ contains $y^2$. We may then define the *generality* of a polynomial model as its *highest polynomial variable,* which are $x^2$ for $M_1$ and $y^2$ for $M_2$.

If the interpretations of x and y are fixed and cannot be swapped, then we cannot say that one is more general than the other and we can call them incomparable in generality.

The described measures are computed separately for each individual model rather than for a pair of models to be compared. As a result, the measures may not represent an *intuitive* generality order relation between models. For instance, we can call model $M_3$ *more general* than model $M_2$, $M_3 >_{Mg} M_2$, because $M_2$ is a *specialization* of $M_3$ with b=4. Similarly, intuitively Pmodel

$$M_4: \forall x,y \ ax+cx^2+by^2=0$$

is more general than Pmodels $M_1$, $M_2$, and $M_3$, because coefficients in $M_1$, $M_2$ and $M_3$ are different numeric specializations of a, b, and c in $M_4$.

However,

$$HP(M_1)= HP(M_2)= HP(M_3)=HP(M_4)=2;$$

that is, all of these Pmodels have the same HP generality, while $M_4$ is intuitively more general than the other models.

Thus, alternatively, we may define the *generality* of a polynomial model as its *highest polynomial variable,* which is $x^2$ for $M_1$ and $y^2$ for $M_2$. If the interpretations of x and y are fixed and we cannot swap symbols x and y, then we cannot say that one of them is more general than the other. Thus, they would be incomparable in generality.

Definition. *HPV measure of polynomial model generality* is defined as the *Highest Power Variable (HPV)* of the polynomial model.

For models $M_1$ and $M_2$ we have $HPV(M_1)=x^2$ and $HPV(M_2)=y^2$. If there are two HPV as in model $M=2x^2+4x+3y^2+5y=0$, then HPV(M) is a pair $(x^2,y^2)$. In the case with more than two HPVs we will have an n-dimensional vector of HPVs.

Below we discuss the advantages and disadvantages of HP and HPV measures. If we cannot fix the meaning of the variables x and y, or if we cannot even agree to use only the symbols x and y, then HPV may not be an appropriate measure.

Swapping symbols will lead to $HPV(M_1)=y^2$ and $HPV(M_2)=x^2$. Using other symbols such as z and v instead of x and y will make a use of HPV even more questionable.

Both HP and HPV are measures that are applied to individual models, not to pairs of models. As a result, they may likewise not represent an *intuitive* generality relation between Pmodels. For instance, we can call $M_3$ *more general* than $M_2$, $M_3 >_{Mg} M_2$, because $M_2$ is a specialization of $M_3$ with b=4. Similarly, intuitively model

$$M_4: \forall x,y \; ax+cx^2+by^2=0$$

is more general than models $M_1$, $M_2$, and $M_3$, because the coefficients in $M_1$, $M_2$ and $M_3$ are different numeric specializations of a, b and c in $M_4$.

The HPV measure captures this better with

$HPV(M_1)=x^2$, $HP(M_2)=y^2$, $HP(M_3)=y^2$, $HP(M_4)=(x^2,y^2)$.

$M_4$ is more general then other models in HPV, but model $M_3$ is a more general than $M_2$, which is a special case of $M_3$ with b=4.

This is not captured by HPV,

$$HPV(M_2)=HPV(M_3)=y^2.$$

Therefore, we introduce another generality characteristic that is defined on pairs of models.

<u>Definition</u>. Polynomial model $M_i$ is a coefficient specialization (*C-specialization*) of a polynomial model $M_j$ if coefficients of $M_i$ are specializations of coefficients of $M_j$. In other words, $M_j$ is a *C-generalization* of model $M_i$.

For instance, Pmodel $M_3$ is a C-specialization of model $M_4$. Similarly, $M_2$ is a C-specialization of $M_3$, but $M_1$ is not C-specialization of $M_2$.

Note that model $M_4$ is more uncertain than models $M_1$, $M_2$, and $M_3$, because all coefficients in $M_4$ are uncertain, but none of the coefficients are uncertain in $M_1$, $M_2$, and only one coefficient is uncertain in $M_3$.

<u>Definition</u>. SP *measure of polynomial model simplicity* is defined as the sum of powers of its variables. For instance, $SP(M_3)=3$ having powers 1 and 2 in $5x^1 + by^2$. One model is SP-simpler than the other if its SP measure is smaller. Here $M_3$ is simpler than $M_4$, which has SP=5.

Uncertainty, generality, and simplicity relations can be *isomorphic* (produce the same order of Pmodels), or be quite different. In the next section, we provide a parameterization mechanism that highlights the difference.

## IV. PARAMETERIZATION

Below we parameterize uncertainty and generality of polynomial models.

The first block in Table 1 shows a chain of models with increasing NUC uncertainty, (from 0 to 2), fixed PH generality, (2), fixed HPV generality, ($y^2$), and increasing C-generality: that is, from all known coefficients to two unknown coefficients, a and b, that substitute coefficients 3 and 4, respectively.

All these models have HP level 2 and SP simplicity level equal to 4. Other blocks in Table 1 illustrate other relations between these characteristics of the models.

Consider another example of the five models at five NUC increasing uncertainty levels, from 0 to 4:

$M_0: \forall x,y \; 9x^2+3y+7x =10$
$M_1: \forall x,y \; x^2+3y+7x =10$
$M_2: \forall x,y \; ax^2+3y+7x = 10$
$M_3: \forall x,y \; ax^2+by+7x = 10$
$M_4: \forall x,y \; ax^2+by+cx+d = 10$

Table 1

| | |
|---|---|
| 1. Chain of models with increasing both NUC uncertainty and C-generality, with fixed HP and HPV generalities and fixed SP simplicity. | $\forall x,y \; 3x+4y+5y^2=0$ <br> $\forall x,y \; ax+4y+5y^2=0$ <br> $\forall x,y \; ax+by+5y^2=0$ |
| 2. Chain of models with increasing NUC uncertainty, without increasing C-generality, and with fixed HP and HPV generalities. | $\forall x,y \; 3x+4y+5y^2=0$ <br> $\forall x,y \; ax+9y+5y^2=0$ <br> $\forall x,y \; ax+by+7y^2=0$ |
| 3. Chain of models with increasing NUC uncertainty, without increasing C-generality, and increasing HP generality, but with incomparable HPV generality. | $\forall x,y \; 3x+4y+5y^2=0$ <br> $\forall x,y \; 3x^2+by \;\;\;=0$ <br> $\forall x,y \; ax+7y+by^2=0$ |
| 4. Chain of models with increasing NUC uncertainty and SP simplicity, but with decreasing HP generality, incomparable HPV generality and without increasing C-generality (ax does not generalize $3x^2$). | $\forall x,y \; 3x+4y+5y^2=0$ <br> $\forall x,y \; 3x^2+by \;\;\;=0$ <br> $\forall x,y \; ax \;+by \;\;\;=0$ |

Models $M_4$, $M_3$, $M_2$, $M_1$, $M_0$ form an *uncertainty decreasing chain* on the UNC uncertainty relation defined above: $M_4 >_{Mu} M_3 >_{Mu} M_2 >_{Mu} M_1 >_{Mu} M_0$. These models also form a *generality decreasing chain* $M_4 >_{Mg} M_3 >_{Mg} M_2 >_{Mg} M_1 >_{Mg} M_0$. Here every model $M_i$ can be obtained by specialization of parameters of model $M_{i+1}$.

Each model has four parameters, $p_1$, $p_2$, $p_3$, and $p_4$. For instance, for model $M_2$ parameter $p_1=1$ represents uncertainty of $ax^2$, with an unknown coefficient a. Similarly, $p_2=p_3=p_4=0$, because further coefficients 3, 7 and 10 are known.

Thus, each model is parameterized as a Boolean vector, $\mathbf{v_i} =(v_{i1},v_{i2},..,v_{ik},…, v_{in})$: $M_4$: $v_4=1111$; $M_3$: $v_3=1110$; $M_2$: $v_2=1100$; $M_1$: $v_1=1000$; $M_0$: $v_0=0000$.

<u>Definition</u>. Parametric model $M_i$ is *no less general* than model $M_j$ if $\mathbf{v_i} \geq \mathbf{v_j}$, i.e., $\forall k \; v_{ik} \geq v_{jk}$.

In accordance with this definition, we have

$$1111 \geq 1110 \geq 1100 \geq 1000 \geq 0000,$$

which is isomorphic to

$$M_4 \geq_{Mg} M_3 \geq_{Mg} M_2 \geq_{Mg} M_1 \geq_{Mg} M_0.$$

In other words, the model with parameters (0000) is a specialization of all previous models.

**Learning Operator**. The intuitive idea of learning Pmodels from data is to get a more specialized model from a given Pmodel M.

<u>Definition</u>. Mapping C is called a **learning (adaptation) operator** C,

$$C(M_i,E)= M_j.$$

This operation represents a *cognitive learning process* C of a new model $M_j$ from a given model $M_i$ and data E. In other words, the process is an adaptation of model $M_i$ to data E to produce model $M_j$.

A learning operator $C(M_i,E)$ is applied multiple times to produce a chain of models, where each subsequent model is more specific than the previous model.

$$M_i >_{Mg} M_{i+1} >_{Mg} M_{i+2} \ldots M_{i+k-1} >_{Mg} M_{i+k}$$

Theoretically here we may have "≥" relation instead of ">", but non-trivial learning operator should produce more specific models, e.g., we might start from the model with all parameters uncertain (111111) and end up with a more certain model parameterized with a Boolean vector (101010):

$$111111 > 011111 > 011110 > 101010.$$

<u>Definition</u>. Parametric model $M_i$ is *no less general* than model $M_j$ if $\mathbf{u_i} \geq \mathbf{u_j}$, i.e., $\forall k\; u_{ik} \geq u_{jk}$.

Above we encoded known parameters as 1 and unknown as 0. A more detailed uncertainty parameterization is when Boolean vectors are substituted by k-valued vectors $\mathbf{u_i} = (u_{i1}, u_{i2}, \ldots, u_{ik}, \ldots, u_{in+m})$ with

$$u_{ij} \in U = \{0, 1/(k-1), 2/(k-1), \ldots k-2/(k-1), 1\}.$$

V. SIMILARITY MAXIMIZATION

**A similarity maximization** problem is a major mechanism of DLP that is formalized below.

<u>Definition</u>. A similarity $L_{fin}$ measure is called a **final similarity measure** if:

$$\forall\; M, E, L_i\quad L_i(M, E) \geq_{Lu} L_{fin}(M, E).$$

The final similarity measure specifies the level of certainty of model similarity to the data that we want to reach.

<u>Definition</u>. The **static model optimization problem (SMOP)** is to find a model $M_a$ such that

$$L_{fin}(M_a, E) = \max\nolimits_{j \in J} L_{fin}(M_j, E) \qquad (1)$$

subject to conditions (2) and (3):

$$\forall\; M_j \in U(M_a)\; L_{fin}(M_a, E) = L_{fin}(M_j, E) \Rightarrow M_a \geq_{Mu} M_j, \qquad (2)$$

$$\forall\; M_j \in G(M_a)\; L_{fin}(M_a, E) = L_{fin}(M_j, E) \Rightarrow M_a \geq_{Mg} M_j \qquad (3)$$

The goal of (2) and (3) is to prevent model overfitting with data E. Sets $U(M_a)$ and $G(M_a)$ contain Pmodels that are comparable with $M_a$ relative to uncertainty and generality, respectively. Condition (2) means that if $M_a$ and $M_j$ have the same similarity measure with E, then uncertainty of $M_a$ should be no less than uncertainty of $M_j$. Condition (3) expresses an analogous condition for the generality relation. Too specific models can lead to overfitting.

<u>Definition</u>. **The dynamic logic model optimization (DLPO)** problem is to find a Pmodel $M_a$ such that

$$L_a(M_a, E) = \max\nolimits_{j \in J} L_j(M_j, E) \qquad (4)$$

subject to conditions (5) and (6):

$$\forall\; M_j \in U(M_a)\; L_a(M_a, E) = L_j(M_j, E) \Rightarrow M_a \geq_{Mu} M_j, \qquad (5)$$

$$\forall\; M_j \in G(M_a)\; L_j(M_a, E) = L_j(M_j, E) \Rightarrow M_a \geq_{Mg} M_j \qquad (6)$$

This is a non-standard optimization problem. In standard optimization problems, only models $M_i$ are changed but the optimization criterion L is held fixed, since it does not depend on the model $M_i$. In DLP, however, the criterion L changes dynamically with Pmodels $M_j$.

Since the focus of DLP is cutting **computational complexity** (CC) of model optimization, a **dual optimization problem** can be formulated.

<u>Definition</u>. An optimization problem of finding a **shortest sequence** of matched pairs $(M_i, L_i)$ of Pmodels $M_i$ and optimization criteria (similarity measures) $L_i$ that solves the optimization problem (4)-(6) for the given data E is called a **dual dynamic logic model optimization (DDLMO)** problem, which finds a sequence of n matching pairs

$$(M_1, L_1),\; (M_2, L_2), \ldots, (M_n, L_n),$$

such that

$$L_n(M_n, E) = \max\nolimits_{i \in I} L_i(M_i, E)$$

and

$$\forall\; M_i\; L_i = F(M_i),\; C(M_i, E) = M_{i+1},$$

$$M_i \geq_{Mu} M_{i+1},\; M_i \geq_{Mg} M_{i+1},\; M_n = M_a,\; L_n = L_a.$$

This means finding a sequence of more specific and certain Pmodels for the given (1) data E, (2) matching operator F, and (3) learning operator C to maximize similarity measure $L_i(M_i, E)$.

VI. MONOTONE BOOLEAN FUNCTIONS

<u>Definition</u>. A Boolean function $f: \{0,1\}^n \to \{0,1\}$ is a **monotone Boolean function** if: $\mathbf{v_i} \geq \mathbf{v_j} \Rightarrow f(\mathbf{v_j}) \geq f(\mathbf{v_i})$.

This means that $(\mathbf{v_i} \geq \mathbf{v_j}\; \&\; f(\mathbf{v_i}) = 0) \Rightarrow f(\mathbf{v_j}) = 0$ and $(\mathbf{v_i} \geq \mathbf{v_j}\; \&\; f(\mathbf{v_j}) = 1) \Rightarrow f(\mathbf{v_i}) = 1$. Function f is non-decreasing. Consider fixed E and $M_i$ that are parameterized by $\mathbf{v_i}$ and interpret $L(M_i, E)$ as $f(\mathbf{v_i})$, i.e., $L(M_i, E) = f(\mathbf{v_i})$. Assume that $L(M_i, E)$ has only two values (unacceptable 0, and acceptable 1). It can be generalized to a k-value case if needed. If $L(M_i, E)$ is monotone, then $\mathbf{v_i} \geq \mathbf{v_j} \Rightarrow L(M_i, E) \geq L(M_j, E)$, e.g., if $L(M_{3-1110}, E) = L(M_{2-1100}, E) = 0$, then

$$(\mathbf{v_i} \geq \mathbf{v_j}\; \&\; L(M_{3-1110}, E) = 0) \Rightarrow L(M_{2-1100}, E) = 0 \qquad (7)$$

$$(\mathbf{v_i} \geq \mathbf{v_j}\; \&\; L(M_{2-1100}, E) = 1) \Rightarrow L(M_{3-1110}, E) = 1 \qquad (8)$$

This means that if a model with more unknown parameters $\mathbf{v_i}$ failed, then a model with less unknown parameters $\mathbf{v_j}$ will also fail. If we conclude that a quadratic polynomial model $(M_2)$ is not acceptable, $L(M_2, E) = 0$, then a more specific quadratic model $M_3$ also cannot be acceptable, $L(M_3, E) = 0$. Thus, we do not need to test model $M_3$. This monotonicity property helps to decrease computational complexity. We use

$$(\mathbf{v_i} \geq \mathbf{v_j}\; \&\; f(\mathbf{v_i}) = 0) \Rightarrow f(\mathbf{v_j}) = 0 \qquad (9)$$

for rejecting models, and

$$(\mathbf{v_i} \geq \mathbf{v_j}\; \&\; f(\mathbf{v_j}) = 1) \Rightarrow f(\mathbf{v_i}) = 1 \qquad (10)$$

for confirming models. In the case of a **model rejection test** for data E, the focus is not to quickly build a model but rather to quickly reject a model M—in the spirit of Popper's falsification principle. In essence, the test $L_3(M_3, E) = 0$ means that the whole class of models $M_3$ with 3 unknown parameters fails. Testing $M_3$ positively for data E requires finding 4 correct parameters. This may mean searching in a large 4-D parameter space $[-100, +100]^4$ for a single vector, say $(p_1, p_2, p_3, p_4) = (9, 3, 7, 10)$, if each parameter varies in the interval $[-100,100]$. For rejection we may need only 4 training vectors (x,y,u) from data E and 3 other test vectors. The first four vectors allow us to build a quadratic surface in 3-D as a model. We then simply test whether three test vectors from E fail to fit this quadratic surface.

## VII. SEARCH PROCESS

In the optimization process, we want to keep track of model rejections and to dynamically guide what model will be tested next in order to minimize the number of tests. Formulas (9) and (10) are key equations to minimize tests, but we need the whole strategy for how to minimize the number of tests and to formalize it. This strategy is formalized as minimization of **Shannon function** $\varphi$ [1]

$$\min_{A \in \mathbf{A}} \max_{f \in \mathbf{F}} \varphi(f,A),$$

where **A** is a set of algorithms, F is a set of monotone functions, and $\varphi(f,A)$ is a number of tests that algorithm A does to fully restore function f. Each test means computing a value f(**v**) for a particular vector **v**. In the theory of monotone Boolean functions it is assumed that there is an **oracle** that is able to produce the value f(**v**), thus each test is equivalent to a request to the oracle [1,4,5]. Minimization of the Shannon function means that we search for the algorithm that needs the smallest number of tests for its worst case (function f that needs maximum number of tests relative to other functions). This is a classic min-max criterion. It was proved in [1] that

$$\min_{A \in \mathbf{A}} \max_{f \in \mathbf{F}} \varphi(f,A) = \binom{n}{\lfloor n/2 \rfloor} + \binom{n}{\lfloor (n/2) \rfloor +1},$$

where $\lfloor x \rfloor$ is a floor of x. The proof is based on the structure called *Hansel chains*. These chains cover the whole n-dimensional binary cube $\{1,0\}^n$. The steps of the algorithm are presented in detail in [3,5]. The main idea of these steps is building Hansel chains, starting from testing the smallest chains, expanding each tested value using formulas (7) and (8), testing values that are left unexpanded on the same chains, then moving to larger chains until no chains are left. The goal of the search is to find a **smallest lower unit v**, i.e., a Boolean vector such that f(**v**)=1, and for every **w**<**v** f(**w**)=0, and for every **u**>**v** |**u**|>|**v**|. A simpler problem could be to find any lower unit of f.

The search problem in logical terms can be formulated as a satisfiability problem: Find a system of axioms $T_a$ such that $L_a(T_a, E)=1$ subject to the condition

$$\forall\; T_j\;\; L_a(T_a,E) = L_j(T_j, E) \Rightarrow T_j \geq_{Mu} T_a,$$

i.e., if $T_a$ and $T_j$ have the same similarity with E, then $M_a$ should have a lower uncertainty than $T_j$, e.g., $T_j \geq_{Mu} T_a$. If a similarity measure L is defined as a **probabilistic measure** in [0,1] then the probabilistic version of the task of finding a system of axioms T for model A that maximizes a probabilistic similarity measure is:

$$\text{Max}_{i \in I}\; L(M_i, E).$$

## VIII. DYNAMIC LOGIC VISUALIZATION

Visualization of DL helps to monitor, control, and understand the DL search process. It allows us to explore visually the relations between the verified and refuted models, and to discover dynamic patterns of model search. A model to be tested is shown as a highlighted bar (see Figure 2). When the answer (verified or refuted) is provided, this bar changes its color to red (refuted) or black (verified). Next, these bars are expanded using monotonicity automatically or by a user. A user can see simultaneously all models that have been currently tested as bars. A red bar (marked as model $M_1$) is a refuted model $M_1$ and yellow bars under $M_1$ are models refuted using monotonicity. Similarly, model $M_2$ and $M_3$ were tested next and refuted, which is indicated by their red bar. Colored bars below them are also refuted using monotonicity. Then model $M_4$ was verified and encoded by a black bar with all models verified using monotone expansion of $M_4$ are shown above $M_4$ as grey bars. Using technique from [Kovalerchuk, Delizy, 2005], bars can be rearranged to reveal their pattern as a Pareto set/border (see figure 2).

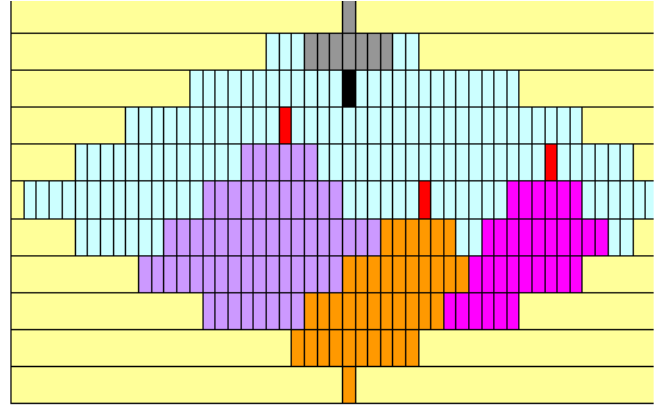

Figure 2. Dynamic logic process visualization

The border is dynamically developed as more models are tested. Finally, a user will see a border between verified and refuted models and can impose additional requirements to monitor border models that satisfy such additional requirements.

## IX. FORMAL DESCRIPTION OF P-DYNAMIC LOGIC

Below we describe formally our dynamic logic of models of phenomena (DLP) as a summarization of concepts introduced above. As with any logic, this logic consists of three parts: (1) a semantic part, (2) a syntactic part, and (3) a reasoning part.

### A. Semantic part of DLP

At first, we define the **semantic part** of DLP. It consists of two related algebraic systems. The first one is as a *three-sort*

*algebraic system* [7], where {E} are sets of data, {M} are sets of Pmodels, R is the set of real numbers,

$$\mathbf{EM} = <\{E\}, \{M\}, R; \Omega_{EM}>. \quad (11)$$

$\Omega_{EM}$ consists of sets of relations, $\Omega_E$, $\Omega_{PM}$, $\Omega_R$ on {E}, {M} and R, respectively, and operators {L} and C that connect {E}, {M} and R,

$$\Omega_{EM} = <\Omega_E, \Omega_{PM}, \Omega_R, \{L\}, C>. \quad (12)$$

Here $\Omega_{PM} = \{\geq_{Mg}, \geq_{Mu}, \geq_{Ms}\}$ presents partial order relations between Pmodels relative to their generality, uncertainty and simplicity as described above. Each **similarity (correspondence) measure** $L_i \in \{L\}$ is a mapping:

$$L_i: \{M\} \times \{E\} \to R \quad (13)$$

that captures numerically the similarity between data and P-models. A higher value of $L_i(M,E)$ indicates a higher consistency between data E and the Pmodel M in the aspect that was captured by $L_i$. Having a set of such measures {L} allows us to choose appropriate measures dynamically.

Next, a Pmodel **enhancement (learning) operator C** is

$$C: \{M, E\} \to \{M\} \quad (14)$$

that changes a P-model. This operator brings us *dynamics* of model changes. Thus, we have two types of dynamics and we need a formal way to express the change of L similar to (14) for models. This is done by introducing a *two-sort algebraic system*

$$\mathbf{ML} = <\{M\}, \{L\}; \Omega_{PM}, \Omega_L, F>, \quad (15)$$

with

$$\Omega_{PM} = \{\geq_{Mg}, \geq_{Mu}, \geq_{Ms}\}, \quad \Omega_L = \{\geq_{Lg}, \geq_{Lu}\},$$

and F as a mapping between sets {M} and {L} that preserves relations on {M},

$$F: \{M\} \to \{L\}. \quad (16)$$

Algebraic systems (11) and (15) have several common components, but {L} has very different roles in each. In (11), {L} is part of Ω, i.e., it is an *operator*; but in (15) {L} is one of two *base sets* of the algebraic system, where the operator is F: {M} → {L}.

Thus, we cannot simply join systems (11) and (15) together to a single algebraic system because of the different roles played by {L} in (11) and (15). We would need a generalized multi-sort algebraic system in second order-logic to do this. Note that separately (11) and (15) are much simpler systems in first-order logic (FOL).

In addition, we can build (15) without a specific dataset E∈{E}, because (15) does not require {E}. To construct (11) for a specific dataset E, we would need to define $\Omega_E$. See an extensive discussion of these issues and examples in [6,4]. Thus, the semantic part of the DLP is

$$\mathbf{EMML} = <\mathbf{EM}, \mathbf{ML}; W>,$$

where $W(M_i, E, M_j)$ is the following relation that consists of three parts (17)-(19):

$$W(M_i, E, M_j) \equiv$$
$$([C \in \Omega_{EM} \& C(M_i,E) = M_j] \& \quad (17)$$
$$\{[(M_i >_u M_j) \vee (M_i >_g M_j) \vee (M_j >_s M_i)] \vee \quad (18)$$
$$[L_i = F(M_i) \& L_j = F(M_j) \& L_j(M_j,E) > L_i(M_i,E)]\}) \quad (19)$$

Relation $W(M_i, E, M_j)$ is true if and only if C produces a Pmodel $M_j$ using E that is better than the input P-model $M_i$ in at least one of its characteristics (i.e., more certain, more specific, simpler, or better fit data relative to similarity measures). In other words, a Boolean predicate $W(M_i, E, M_j)=1$ if $C(M_i,E)=M_j$ produces an improved model $M_j$, otherwise $W(M_i, E, M_j)=0$.

In accordance with the definitions in previous sections, in (18) $M_i >_{Mu} M_j$ means that $M_j$ is a more certain model than $M_i$. Similarly, $M_i >_{Mg} M_j$ means that $M_j$ is a more specific, and $M_j >_s M_i$ means that $M_j$ is simpler than $M_i$. Property (19) means that model $M_j$ better fits data E than model $M_i$ relative to measures $L_j$ and $L_i$, which are dynamically assigned to $M_j$ and $M_i$ by applying F to them.

Not every operator that produces another model can be called a learning operator. Relation W sets up a semantic criterion for operator C to be a learning operator for P-models $M_i$, $M_j$ and data E. Note that (17) and (18) in W can be checked having only **EM**, but (19) requires **ML** too.

Now we want to discuss how to compute truth-values of predicates and outputs of operators in **EM** and **ML**. For **EM** and **ML** we gave some examples of $\Omega_M$ and $\Omega_L$, {L}, C and F. An extensive set of examples of similarity measures L and learning operators C are given in [8]. For ML in addition we need to compute F. For instance, if M(c,r) is a model that represents a circle with center c and radius r, then L(c,r) could be the same, F(M(c,r)) = L=M(c,r). This means that L tests exactly model M(c, r). An alternative F could be F(M(c,r))= L=M(c, r+e) that accepts all models with radiuses no greater than r+e.

In general a library of such matching operators F and relations $\Omega_M$ and $\Omega_L$ should be created that will be available to researchers. Thus, we have a semantic part of the P-dynamic logic identified.

*B. Syntactic part of P-Dynamic Logic*

Now we need to explore which part of this semantic machinery can be transferred to the syntactic level so that we can do reasoning without semantic knowledge to get at least some non-trivial inferences and conclusions.

To distinguish the syntactic level from the semantic level, we will use low-case notation **em** and **ml** to define the syntactic structure of **EM** and **ML**, respectively. We do this for all components; thus,

$$\mathbf{em} = <\{e\}, \{m\}, r; \omega_{em}>. \quad (20)$$

The $\omega_{em}$ consists of sets of relations $\omega_e$, $\omega_m$, $\omega_r$ on {e}, {m} and r, respectively, and operators {l} and c that connect {e}, {m} and r,

$$\omega_{em} = <\omega_e, \omega_m, \omega_r, \{l\}, c> \quad (21)$$

Similarly, we define

$$\mathbf{ml} = <\{m\}, \{l\}; \omega_m, \omega_l, f>, \quad (22)$$

that includes

$$\omega_m = \{\geq_{mg}, \geq_{mu}, \geq_{ms}\}, \omega_l = \{\geq_{lg}, \geq_{lu}\} \; f: \{m\} \to \{l\} \quad (23)$$

Similarly, we define $\mathbf{emml} = <\mathbf{em}, \mathbf{ml}; w>$, where

$$w(m_i, e, m_j) \equiv$$
$$( [c \in \omega_{em} \;\&\; c(m_i, e) = m_j] \Leftrightarrow \quad (24)$$
$$\{[(m_i >_{mu} m_j) \vee (m_i >_{mg} m_j) \vee (m_j >_{ms} m_i)] \vee \quad (25)$$
$$[l_i = f(m_i) \;\&\; l_j = f(m_j) \;\&\; l_j(m_j, e) > l_i(m_i, e)]\}) \quad (26)$$

*C. Reasoning part of DLP*

Now we will discuss the reasoning part of DLP. To build a reasoning *deductive system* for DL we need a set $\Lambda$ of *logical axioms*, a set $\Sigma$ of *non-logical axioms*, and a set $\{\Gamma, \varphi\}$ of *rules of inference*.

It will be interesting to explore if it is *complete*, i.e., for all formulas $\varphi$,

If $\Sigma \; \Box \; \varphi$ then $\Sigma \; \Box \; \varphi$.

If this system is actually incomplete, then we will have statement $\varphi$ such that neither $\varphi$ nor $\neg\varphi$ can be proved from the given set of axioms. In contrast, in a complete system, all true statements (made true by the set of axioms) are provable.

We also need to know the *validity* of the inference rules that is, that the conclusion follows the premises. The soundness of these rules needs to be clarified too; that is, it needs to be shown that the conclusion follows from the premises when the premises are in fact true.

While all these are interesting research issues, at this stage of formalization of DLP it is more important to have a rich set of useful non-logical axioms $\Sigma$ first.

To clarify this issue we need to answer the following questions: "What reasoning is possible in DLP?" What could be the most interesting part of syntactic reasoning for a specific practical problem?"

Assume that we already have a knowledge base (KB). This KB contains $m_1,\ldots m_n$, and e that are called facts in this KB ($m_1,\ldots m_n$ are interpreted semantically as Pmodels and e as data). KB also contains some expressions, e.g., $w(m_i, e, m_j)$ (interpreted semantically as Pmodel $m_j$ is an improved Pmodel $m_i$).

For the first series of questions we have common first-order reasoning with logical axioms $\Lambda$ that use $\wedge, \vee, \neg$, and $\Rightarrow$ operators.

Non-logical axioms $\Sigma$ include the disjunction axiom (DA) and conjunction axiom (CA) and mixed CDA axiom:

CA: $w(m_i, e_1, m_k) \vee w(m_i, e_2, m_k) \Rightarrow w(m_i, e_1 \cup e_2, m_k)$,

DA: $w(m_i, e_1, m_k) \wedge w(m_i, e_2, m_k) \Rightarrow w(m_i, e_1 \cap e_2, m_k)$,

DCA: $w(m_i, e_1, m_k) \wedge w(m_i, e_2, m_k) \Rightarrow w(m_i, e_1 \cup e_2, m_k)$,

The last axiom is redundant if we add the inclusion axiom (IA):

IA: $w(m_i, e_1 \cap e_2, m_k) \Rightarrow w(m_i, e_1 \cup e_2, m_k)$

Together DA and IA produce DCA. The real world interpretation of these axioms will assume some regularity in data and learning operators. They should be "reasonable."

Next we want to know if $w(m_i, e, m_k)$ can be inferred purely syntactically knowing that $w(m_i, e, m_j)$ and $w(m_j, e, m_k)$ are in the KB but without using a semantic interpretation of $w(m_i, e, m_k)$. Having this we will have useful syntactic reasoning in this DL. In fact, such *transitivity axiom* (TA) can be established as a part of a set of non-logical axioms $\Sigma$,

TA: $w(m_i, e, m_j) \wedge w(m_j, e, m_k) \Rightarrow w(m_i, e, m_k)$

It follows from the semantics of relation W that it is transitive, thus we can postulate transitivity for w. As a result of this postulate we can infer $w(m_i, e, m_k)$ purely syntactically having $w(m_i, e, m_j)$ and $w(m_j, e, m_k)$ in KB without going to the semantic level and computing $W(M_i, E, M_k)$, which can be computationally challenging for a large dataset E. This is a major advantage of using DLP syntactic reasoning instead of computations at the semantic level.

The reasoning mechanism $T_{em}$ in **em** is a first-order logic with terms in its signature. Similarly, $T_{ml}$ in **ml** it is a first-order logic in terms of its signature. In $<\mathbf{em}, \mathbf{ml}; w>$ we have $T_{emml}$ that is also a first-order logic reasoning with w, but if we substitute w with its components (24)-(26), we will have second-order logic reasoning. Thus, the complete description of the **dynamic logic of phenomena models** is

$$DL = <DL_{EM}, DL_{ML}, DL_{EMML}>,$$

where

$$DL_{EM} = <\mathbf{em}, T_{em}, \mathbf{EM}>,$$
$$DL_{ML} = <\mathbf{ml}, T_{ml}, \mathbf{ML}>,$$
$$DL_{EMML} = <\mathbf{emml}, T_{emml}, \mathbf{EMML}>.$$

*D. Learning operator*

Now we want to elaborate the concept of a learning operator, C. There is an important question about this operator: "How sophisticated should C be relative to a brute force algorithm that computes L(M,E) for every P-model M and selects M that provides min L(M,E)?"

It will be very advantageous to get a simple and quite universal operator, C, which will allow us to use it for solving a wide variety of problems. Let us consider a direct modification of brute force algorithm to the situation with dynamic change of correspondence (similarity) measures $L_i$, which is assumed in the dynamic logic. The next assumption is that the space of highly uncertain models is relatively small. Thus, a brute force algorithm can work for these models in a reasonable time. Next, the best model $M_1$ found at this step will produce a new correspondence measure $L_1$ and this measure will be applied to a *new set of models* produced by $M_1$ for evaluation (e.g., with a more dense grid around $M_1$ or in another location if $L(M_1)$ is low). To produce new models we introduce a new operator, H, that we will call a **specialization operator**, as follows:

Step 1. Select initial Pmodel $M_0$.
Step 2. Produce set of Pmodels $H(M_0)$.

Step 3. Compute $L_0(M, E)$ for every M from $H(M_0) = \{M\}$ and find model $M_1 = \arg \min_{\{M\}} L_0(M, E)$.

Step 4. Test if $L_0(M_1, E) > T$, i.e., is above the needed correspondence threshold. Stop if it is true, else go to step 5.

Step 5. Repeat steps 1-4 until all models tested or time limit reached.

More complex strategies for the H operator can be based on breadth-first, depth-first, and branch-and-bound strategies.

## X. RELATIONS WITH OTHER DYNAMIC LOGICS

A link with the Dynamic Logic for Belief Revision [15] is feasible to explore. The transition from model $M_i$ to model $M_j$ using partial orders between models can be viewed as a belief revision.

At a very general level, the link between Dynamic Logic of Phenomena and either Propositional or First-Order Dynamic Logic [26, 27] follows from the fact that Dynamic Logic of Phenomena contains static components (propositional/first-order formulas), and dynamic components (actions/programs). A formula can change when actions are applied, and actions can be repeated and applied consecutively. Thus, the Dynamic Logic of Phenomena can be viewed as a special case of Propositional or First-Order Dynamic Logics with special types of "programs" described in this paper.

Below we follow [13] in summarizing the Propositional Dynamic Logic that has been useful for capturing important properties of programs, such as correctness. The basic concept of this logic is the state S of a system. Each state S has its properties expressed as a propositional formula $\varphi$. These properties are changed when a program (action) $\pi$ is applied to S with a given precondition a. Thus, program $\pi$ can produce S with property $\psi$ starting with S with property $\varphi$ and precondition a,

$$\pi(S, a) = \psi.$$

Propositional formulas are static components of this logic and actions/ programs are dynamic components. Actions can be applied consecutively and repeated, which is expressed syntactically as $\pi$, $\pi'$ and $\pi^*$, respectively. An action also can be taken randomly from a set of actions, $\pi \cup \pi'$. One of the actions is testing, denoted as $\varphi$? It tests if property $\varphi$ is true or false for system S.

There is also an order relation ( $\geq_w$ ) on the set of preconditions $\{a\}$. It expresses the relative strength (weakness) of the precondition. The statement $a_1 \geq_w a_2$ means that $a_1$ is weaker than $a_2$. $[\pi]\varphi$ is the *weakest precondition* $\psi$ that leads to postcondition $\varphi$ after performing $\pi$ on S. Also $[\varphi\ ?]\psi$ denotes the classical implication $\psi \to \varphi$, that is if S satisfies $\psi$ then it satisfied $\varphi$, too.

In DLP, the concept similar to the concept of system S in the Propositional Dynamic Logic is a Pmodel M and data E pair, <M, E>. Actions are operations that allow us to produce a new Pmodel from M and E,

$$\pi(M_i, E) = M_{i+1}$$

The concept of precondition is not defined in DLP, but data E or properties of E can be interpreted as a precondition. An alternative way to specify preconditions is to introduce an additional concept of precondition to DLP and to interpret it for specific tasks.

In DLP we have a sequence of actions $\pi_1, \pi_2, \ldots, \pi_n$ that can be iterations $\pi^*$ of a fixed learning algorithm/program $\pi$. In this iteration, the last model is the "best" one for data D. If the desirable property $\psi$ of the model $M_i$ is defined, then we can find the *minimal (weakest) precondition* $[\pi]\psi$ that is needed to convert the previous model $M_{i-1}$ to $M_i$. When we interpret a precondition as data E, then $[\pi]\psi$ will mean the minimal dataset that is needed to get a model M with property $\psi$. Thus, designing actions $\pi_k$ by using a learning operator $C(M_{i-1}, E)$ and then testing that the result of this action/program $\varphi$ implies $\psi$, $[\varphi\ ?]\psi$, provides a link between DLP and propositional dynamic logic.

In the Logic of Action [16, 28-30], models are transformed based on the information communicated that uses Kripke's model of actions and Kripke models of states [30]. The update of an epistemic model is done by an epistemic action (i.e., an action that affects the epistemic state of a group of agents). An epistemic model can be matched with a phenomena model in DLP, and an epistemic action will be a program that implements a learning operator $C(M, E)$.

The focus of epistemic logic is expressing the distributed, localized information that is accessible only in parts of the system. In Epistemic logic [31,32] we have $\varphi!$ as an action of learning the truth of proposition $\varphi$ that is the result of the truthful *public announcement* of $\varphi$. In DLP, we have a learning operator $C(M, E)$ that learns a new phenomena model that is more complex than a single proposition. Thus, the epistemic logic can be expanded for the phenomena models. Next, the actual process of getting $\varphi!$ can be elaborated in both DLP and the epistemic logic.

Public announcements are presumed to be certain, but often learning a parameter is uncertain. Qualitative logics for "risky knowledge" [33,34] have been proposed, which capture rational but *uncertain acceptance of a proposition*. This approach, based on $\varepsilon$-acceptability [38], combines both probability and modality to characterize tentative or corrigible acceptance of a set of sentences. A sentence $\varphi$ is $\varepsilon$-accepted if the probability of $\neg\varphi$ is at most $\varepsilon$, where $\varepsilon$ is taken to be a fixed small parameter.

Logics of rational acceptance are inspired by classical statistical methods, such as Fisher's exact test for detecting non-random association between variables. Fisher's exact test specifies conditions under which a null, no-effect hypothesis is rejected, which is nominally equivalent to accepting that the association is non-random.

The probabilities associated with this approach are evidential probabilities, developed in [38] and [39], which construes probability as a metalinguistic relation between a set of sentences and a single sentence on analogy to provability. Evidential probability is interval-valued, defined in a first-order language with the capacity for expressing known statistical frequencies. For example, the language can express that a test of a hypothesis of size α that yields a point in the rejection region supports the denial of the null hypothesis $H_0$ to degree [1-α, 1], or runs a risk error at most α

Normal modal logics are inappropriate for qualitative representations of ε-acceptance, since all instances of the schema

$$(\mathbf{C}) \quad ([]A \wedge []B) \rightarrow [](A \wedge B)$$

are valid in normal Kripke models. But if the box modality [] is interpreted as 'has high probability', distributing 'has high probability' across conjunction should not be a valid principle. That event A has probability greater than ε and event B has probability greater than ε does not entail that the joint event of A and B has probability greater than ε.

Fortunately, Schema C is not valid within the minimal models of classical modal logic [35], and belief revision operators can be defined for monotone classical systems without C [36], as well as rudimentary measures of robustness [37].

A connection between Dynamic Epistemic Logic, expanded with ε-acceptability, and DLP can be expressed through a **similarity (or correspondence) measure**, L: {(M,E)} → R, which serves a role similar to probability in Epistemic logic with ε-acceptability.

XI. CONSEQUENCES VERSUS CONDITIONALS

As we observed in the introduction, DLP locates dynamics in the *meta-language* rather than in the *object language*, which is in contrast to Dynamic Epistemic Logic. We considered ways to bridge these two approaches in the last section, closing with a description of a qualitative modeling of ε-accepted "risky knowledge." In that case the modal techniques of revision, common to dynamic epistemic logic, were harnessed for a specific purpose, to evaluate the robustness of ε-accepted sentences. Here we observe reasons for respecting the distinction.

The object-language/meta-language distinction is important to observe [40], since there is a difference between procedural rules and declarative sentences that often needs to be preserved rather than ameliorated. For while declarative sentences take truth values and can be embedded in complex formulas, procedural rules often do not. However, the question of whether an "if…then…" construction is better viewed as a conditional formula than as a conditional rule, or whether there is little difference between the two and an analogue of the deduction theorem can be used, is not always clear cut, and will depend on the purpose of the formalization.

To illustrate, input/output logic [41, 42] was conceived as a response to deontic logics handling of norms and declarative statements. *Norms*, unlike statements, may be respected or flouted, and may be judged from the standpoint of other norms but are not typically evaluated as 'true' or 'false'. Input/output logic conceives of a norm as an ordered pair of formulas, and a *normative system* is a set G of norms. The task for an input/output logic, then, is to prepare information to be passed into G, and to unpack the resulting output. So, abstractly, the set G is a transformation devise for information, and we may characterize an 'output' operator 'Out' by logical properties typical of consequence operators. A formula x is a 'simple-minded output' of G in *context* a, written x ∈ Out(G,a), if there is a set of norms ($a_1$, $x_1$),…,($a_n$,$x_n$) in G such that each $a_i$ ∈ Cn(a) and x ∈ Cn($x_1$ & … & $x_n$), where Cn(a)={$a_i$| a ⊨ $a_i$} is the classical semantic consequence set of a that is a set of all contexts $a_i$ such that every model of a is also a model of $a_i$ [41].

Out(G,a) satisfies three rules: writing (a,x) for x ∈ Out(G,a), they are *strengthening input* (SI), *conjoining output* (AND), and *weakening output* (WO):

From (a,x) to (b,x), whenever a ∈ Cn(b).  (SI)
From (a,x), (a,y) to (a, x∧y).   (AND)
From (a,x) to (a,y), whenever y ∈ Cn(x).  (WO)

Strengthening input means that if context a leads to output x then a stronger context b also leads to x. Similarly, weakening output means that if context a leads to output x then context a also leads to a weaker output y, x ⊨ y, that is y∈Cn(x).

Here the idea is that, while classical logic may be used to 'process' the input and to 'unpack' the output, the operator Out itself does not have either an associated language or a proof theory. Even so, Out may enjoy structural properties commonly attributed to consequence operators. Simple-minded output is the most general characterization of Out, but stronger input/output operators have been studied [41], including a characterization of Poole's default system [43] and a system similar to Reiter's default logic [25].

**Example** [44] Suppose there are two community norms governing social benefits, one which maintains that poor citizens receive a housing subsidy, and the other that elderly citizens receive health insurance, i.e., G = {(Poor, SubsidizedHousing), (Elderly, HealthInsurance)}. Then, the community should also provide a housing subsidy if no-income implies poor is included in the theory, i.e.,

SubsidizedHousing ∈ Out(G,(NoIncome → Poor) & NoIncome), since

Poor ∈ Cn((NoIncome → Poor) & NoIncome) and
SubsidizedHousing ∈ Cn(SubsidizedHousing).

While one may agree that deontic logic has traditionally mishandled norms, it is not clear that there is anything essential about norms, which would preclude them from being represented by a statement within a language. The point instead is that we must first understand the system we wish to model before defining the language we need, rather than the other way around.

A good illustration of this point is the theory of causal modeling [45,46], since "if…then…" statements have long been thought to express a causal relationship between events. Whereas in normative systems the main issue is the apparent non-truth-functional character of norms, the issues with causal modeling concern the asymmetrical nature of causal relationships and, more importantly, the role that intervention plays in our understanding of causal relationships. Certainly, logics of conditionals have abounded; most attempt to ground conditionals of this kind in the semantic features of the mood, aspect, and tense of natural language conditional statements. And some philosophical theories of causality have attempted to build from these studies in formal semantics.

However, in the case of causal modeling, success came once the order of inquiry was reversed. While both equational and causal models rely upon symmetrical equations to describe a system, causal models in addition impose asymmetries by assuming that each equation corresponds to an independent mechanism, which may be manipulated. The theory of causal Bayesian networks [45,46] supplied the mathematical framework for causal modeling, and it did so by ignoring the task of formulating a calculus and focused instead on getting the right mathematical framework. Only recently has there been a proposal for a calculus, and the clarity of Judea Pearl's causal algebra of 'doing' [44] relies upon this background theory.

## XII. CONCLUSION AND FUTURE WORK

This paper has introduced a formalization of the dynamic logic of phenomena in the terms of first-order logic, logic model theory, and the theory of Monotone Boolean functions. This formalization covers the main idea of DLP -- matching levels of uncertainty of the problem/model to levels of uncertainty of the evaluation criterion, which dramatically decreases computation complexity of finding a model.

The core of the formalization is a partial order on the phenomena models and similarity measures with respect to their uncertainty, generality, and simplicity. These partial orders are represented using a set of Boolean parameters within the theory of monotone Boolean functions and can be visualized in line with [3] to monitor and guide the search in the Boolean parameter space in the DL setting.

Further studies are needed to establish a richer set of non-logic axioms and to explore validity, soundness, and completeness of expanded DLP. Such studies may reveal deeper links with classical logic problems such as decidability, completeness, and consistency.

Further theoretical studies may find deeper links with classical optimization search processes, and may significantly advance them by adding an extra layer of optimization criteria constructed dynamically. This work also gives a new perspective for machine learning [2,4] to develop learning algorithms that can learn evaluation criteria and models simultaneously. We expect that deep links with a known Dynamic Logic of programs [12] motivated by logical analysis of computer programs will also be established.

The proposed formalization creates a framework for developing specific applications and modeling tools. We envision a modeling tool that will consist of sets of models, matching similarity measures, processes for testing them and model learning processes for specific problems in pattern recognition, data mining, optimization, cognitive process modeling, and decision making.

From our viewpoint, the most interesting and useful future DLP research is discovering a mechanism for changing models and similarity measures. Humans have demonstrated capabilities to switch evaluation criteria instantaneously in dynamic environments [11]. This is the area where logic, mathematical modeling, and cognitive science can provide mutual benefits in discovering mechanisms for changing models and changing similarity measures.

APPENDIX 1: COMPUTATIONAL COMPLEXITY

*E. Polynomial complexity case*

Below we describe two examples that give computational motivation to our research on building a logical formalism for representing the dynamic logic of phenomena. The first example is about finding a model to fit a hidden circle in noisy data when the level of noise is so high that a direct human observation does not allow a visual clue about the location of the circle within the image. In the second example, the shape of the object (model) to be found is more complex (see Figure 3 on the right). To identify the circle we need two points, and to identify the second (parabolic) shape we need 4 points A,B,D and G (see Figure 3), that is models with 2 and 4 parameters, respectively.

To solve the first example we assume:
1. The hidden shape to be discovered is a circle.
2. The density of the points in the circle is higher than outside of the circle W, but it is so subtle that it is below of the human perception abilities.
3. The *center* c=($x_c$, $y_c$) of the circle *is not known*, but its range is known: $x_c \in [x_{min}, x_{max}]$, $y_c \in [y_{min}, y_{max}]$.
4. The *resolution* of x and y *are known*, $e_x$ and $e_y$, that is x has a finite set of $n_x$ values, $x_{min}$, $x_{min}+e_x$,…, $x_{min}+ke_x$,…$x_{max}$ and y has a finite set of $n_y$ values, $y_{min}$, $y_{min}+e_y$,…,$y_{min}+re_y$,…$y_{max}$
5. The *radius* R of W *is not known*.
6. The radius R accuracy resolution is $e_r$ =min($e_x$,$e_y$).
7. More than one circle can be in the area and the *number* of the circles *is not known*.
8. Different circles may or may not overlap.
9. The circles' area is less than 1/3 of the whole area.

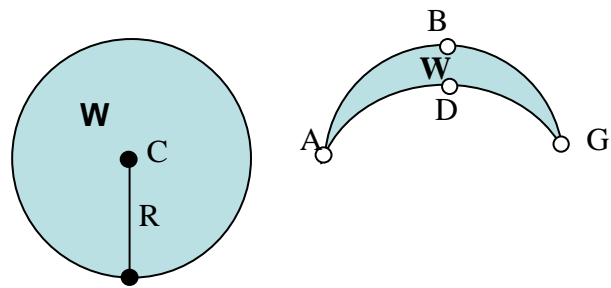

Figure 3. Modeled shapes

The brute force algorithm for example 1 conducts Density Difference Test for every node $C=(x_c,y_c)$ and every R on the grid. It returns 1 if $D(R) - D(\neg R) > T$, else 0, that is if density $D(W)$ in circle W is greater than $D(\neg W)$ with threshold T, where $D(\neg W)$ is density outside W. Here

The **Computational Complexity** (CC) of this algorithm with the base operation as Density Difference Test (C,R) is $O(n^3)$ on a square grid $n \times n$.

The CC of Density Difference Test function is defined by the total number $m$ of given points with the base operation as testing $(x_i - x_c)^2 + (y_i - x_c)^2 \leq R^2$, that is if the point is inside of the circle. The total computational complexity of the brute force algorithm is $O(n^3 m)$. If input points covers the whole grid, then $m=n$ and CC is equal to $O(n^4)$. If fact $m$ is a fraction of $n$, thus if $m=n/10$, then we still have $O(n^4)$.

Now we will change our assumption 1, i.e., that the hidden shape to be discovered is a circle. A more complex shape to discover is produced by two quadratic curves shown in Figure 3. Three points, A, B, and G, are sufficient to identify the first curve, and three points, A, D, and G, are sufficient to identify the second curve. Thus, in total 4 points are sufficient to identify both curves. In the case of the circle, we need only two points, center C, and any point H on the circle. To solve the task for our new parabolic shape, we keep the same assumptions as we had for the circle. We assume the unknown number of shapes of unknown sizes, locations, and orientations. Complexity of a brute force algorithm is $O(n^4)$ which is greater than for Example 1 due to a greater number of parameters involved. The density difference test has the same complexity $O(m)$ as above for the circle, but is more complex for points inside of the shape. However, this test's time grows linearly with $m$. For each point, it computes two quadratic forms and tests if the point is in or out of the shape based on these values. This time is limited by a constant for each given point. Thus, the total complexity is $O(n^4 m)$, and if $m$ would reach $n$, then it will be $O(n^5)$. While these algorithms are polynomial, they have limited applicability for data of practical sizes as we show below.

Assume that $m=n/10$, thus our complexity functions for examples 1 and 2 are $k_1(n^4/10)$ and $k_2(n^5/1)$, respectively. We also assumed $10^9$ base operations/sec in these computations. Note that this base operation is a part of the Density Difference Test computed for each of $m$ inputs.

Tables 1 and 2 show computational complexity under these assumptions for search of circles and parabolic shapes. A computer screen grid of pixels is about 1000×1000 and images from current cameras are larger. Only for a single screen with 100 input points, finding the circle has an acceptable time (100 sec) which also can be too slow for real time applications. Finding a more complex parabolic shape needs 28 hours for the same 100 input points. Increasing $m$ to 1000 or 10000 input points leads to many years as table 2 shows. Real objects, such as aircrafts, have shapes that are more complex. Thus, more critical points are needed to identify their models and CC will be greater. Therefore, fundamentally different approaches are needed to deal with such huge computational complexity issues.

This example illustrates a computational motivation of this paper -- building the dynamic logic formalism for a computationally efficient dynamic logic process.

*E. Exponential complexity case*

Above we considered a situation where we needed to test a polynomial number of models. This was a result of the assumption that models are simple circles or parabolic shapes defined by $n^2$ points on the grid.

Objects with more complex shapes (e.g., long and non-linear tracks) may require identifying more parameters and their combinations that can grow exponentially.

TABLE 1.
COMPUTATION TIME FOR SEARCH OF CIRCLES IN NOISY DATA

| Grid size n | Points m | Circle | |
|---|---|---|---|
| | | Computational Complexity $n^3 m$ | Time for $10^9$ base operations per sec |
| 10 | 1 | 1000 | 0.000001 sec |
| 100 | 10 | 1E+7 | 0.01 sec |
| **1000** | **100** | **1E+11** | **100 sec** |
| **10000** | **1000** | **1E+15** | **277.8 hours** |
| 100000 | 10000 | 1E+19 | 126.8 years |
| 1000000 | 100000 | 1E+23 | 1,27 million years |

TABLE 2.
COMPUTATION TIME FOR SEARCH OF PARABOLIC SHAPES IN NOISY DATA

| Grid size n | Points m | Parabolic shape | |
|---|---|---|---|
| | | Computational Complexity $n^4 m$ | Computational Complexity $n^4 m$ |
| 10 | 1 | 10000 | 0.00001 sec |
| 100 | 10 | 1E+09 | 1 sec |
| **1000** | **100** | **1E+14** | **28 hours** |
| **10000** | **1000** | **1E+19** | **317.1 years** |
| 100000 | 10000 | 1E+24 | 31.7 million years |
| 1000000 | 100000 | 1E+29 | 3.17E+12 years |

APPENDIX 2: EXAMPLE OF PDL PROCESS

Consider a situation with a single unknown object that is described as an interval [a,b] in the total interval [0,10]. We need a process that will find a and b. The model of the object is [a,b], which is highly uncertain because both values a and b are unknown and only limited by the interval [0, 10]. Thus, the first class of models is M={[a,b]: a,b ∈ [0,10]}.

To be specific assume that the unknown object (model) is the interval [0,4]. This model has center c=2 and radius r=2, and can be written as [a,b]=m(c,r). Let also M(c,R)={m(c,r)} be a set of all models with center c and radiuses r that are no greater than R. Consider a set of all models (intervals) M(c,5) with all radiuses r ≤ 5. Here c and r also are highly uncertain.

Now we build at first a very uncertain evaluation measure (criterion) $L_0$ for this class of models M(c,R) and data E as a kernel function. For instance, consider a Gaussian distribution, N(5,10), where c=5 is a mean and r=10 is a standard deviation. This standard deviation is two times greater than the largest R=5 in the [0, 10] interval, thus it covers M(c,5) models (intervals) and should not fail for this class of models. This means that the result of applying $L_0$ based on N(5,10) to M(c,5) is positive, $L_0(M(c,5))=1$. This indicates that an object exists in [0,10], but its location is uncertain. We still do not know exactly c and r, but the class of models M(c,5) is confirmed because we have $L_0(M(c,5))=1$.

Having this positive result, the next step is to change the similarity measure L to make it less uncertain, by substituting $L_0=N(5,10)$ by say $L_1$ based on N(5,7). This can be done by using a learning operator. Testing $L_1$ on M(c,5) is not a computational challenge. It can be done quickly, which is a major benefit of DLP. This quick process of changing L and M is repeated until it reaches the level of maximum certainty of models that is possible on available data. This process is much faster than a brute force approach described in previous examples in this appendix. Here for simplicity of exposition we assumed that L takes binary values. For more complex non-binary L see [8,9].

For the actual use of the DLP methodology, elaborated classes of models and criteria need to be developed at the matched level of uncertainty. Some of them are already developed in likelihood function terms [8,9].